\documentclass[12pt]{article}

\usepackage[margin=1in]{geometry}
\usepackage[T1]{fontenc}
\usepackage[utf8]{inputenc}
\usepackage{lmodern}
\usepackage{microtype}
\usepackage{amsmath,amssymb}
\usepackage{amsthm}
\usepackage{booktabs}
\usepackage{enumitem}
\usepackage{graphicx}
\usepackage[authoryear,round]{natbib}
\usepackage{hyperref}

\hypersetup{
  colorlinks=true,
  linkcolor=blue,
  citecolor=blue,
  urlcolor=blue
}

\newcommand{\bel}{\operatorname{bel}}
\newcommand{\pl}{\operatorname{pl}}
\newcommand{\mpl}{\operatorname{mpl}}
\newcommand{\pr}{\mathsf{P}}
\newcommand{\Unif}{\operatorname{Unif}}

\theoremstyle{definition}
\newtheorem{definition}{Definition}
\theoremstyle{plain}
\newtheorem{theorem}{Theorem}

\title{Inferential Models: The Power of Auxiliary Variables for Reasoning with Scientific Uncertainty}
\author{
Chuanhai Liu\\
\\
Department of Statistics, Purdue University
}
\date{\today}

\begin{document}

\maketitle

\begin{abstract}
A central challenge in scientific inference is to produce uncertainty assessments that are both situation-specific and frequency-calibrated. This article examines inferential models (IMs) as a framework for prior-free probabilistic reasoning with scientific uncertainty. The central IM idea is to view the auxiliary variables in a sampling model as the source of model-based uncertainty. R. A. Fisher's fiducial inference transfers auxiliary randomness to the parameter space before applying probability calculus; IMs instead predict the unobserved auxiliary value with calibrated predictive random sets (PRSs) and transfer the resulting plausibility statements only afterward. This change in order yields valid uncertainty assessments and clarifies the relations among Fisherian fiducial reasoning, Neymanian confidence theory, Dempster-Shafer belief functions, generalized fiducial inference, and IMs. By comparing IMs with objective-prior Bayesian inference, the article argues that E. T. Jaynes' logic-of-science ambition can be continued without forcing all scientific uncertainty into a precise prior distribution because calibrated imprecision is often essential. Finally, the article suggests that a differential-geometric theory of IMs may be within reach, offering a possible route to foundational questions traditionally framed in terms of the likelihood principle.
\end{abstract}


\section{Introduction}

Across its philosophical, methodological, and practical dimensions, as a discipline concerned with transforming empirical experience into knowledge and guiding its application, statistics has long been pulled between two demands. Scientific inference must be conditional on the data actually observed: after seeing the data, scientists seek to learn about their questions of interest from the evidence at hand. At the same time, scientific uncertainty quantification must be frequency-calibrated. A numerical uncertainty statement that is routinely overconfident, or whose scale lacks public meaning, cannot serve reliably as scientific evidence. Bayesian inference satisfies the first demand elegantly when genuine prior information is available. Frequentist inference, in practice, satisfies the second by calibrating procedures through repeated-sampling properties. But these two successes reveal complementary limitations. Bayesian inference requires a joint prior distribution over all unknown quantities, whereas frequentist guarantees are often attached to procedures rather than to direct post-data measures of evidence.

The treatment of Bayes in this review follows from that tension. The aim is not to dismiss Bayesian probability, or Jaynes' broad view in \emph{Probability Theory: The Logic of Science} that probability extends logic to uncertain reasoning \citep{Jaynes2003LogicScience}. It is to separate the formal logic of probability from the scientific interpretation of a particular prior. Information-based priors, such as Jeffreys' invariant prior, Jaynes' maximum entropy priors, and Bernardo's reference priors, are valuable mathematical devices. But when no genuine prior information is available, their posterior probabilities should not automatically be read as calibrated scientific uncertainty. One goal of this article is to give readers enough Inferential Model (IM) background to understand both the power and the limits of such information-based Bayesian prior constructions.

R. A. Fisher's fiducial argument enters at precisely this point \citep{Fisher1935Fiducial,fisher1973statistical}. Fisher rejected Bayesian posteriors with arbitrary priors, but he still wanted post-data uncertainty statements about fixed unknowns. Neyman's confidence theory supplied mathematical clarity by placing probability on the random interval procedure, not on probabilistic statements about the fixed parameter \citep{Neyman1941FiducialConfidence}. Dempster-Shafer theory introduced lower and upper probabilities generated by modeling on the power space of our usual probability space \citep{Shafer1976Evidence,Dempster2008DSCalculus}. Generalized fiducial inference revived the auxiliary-variable inversion in a modern form \citep{Hannig2009GFI}. IMs synthesize these lines by treating the auxiliary variable as the source of uncertainty and by predicting its unobservable realization with calibrated predictive random sets (PRSs) \citep{MartinLiu2013IM,MartinLiu2015Book}.

The present review focuses on IMs and is written from the viewpoint that scientific uncertainty is not just aleatory variation. It is also epistemic: uncertainty about fixed but unknown features of the world, about model adequacy, and about which inferential statements are warranted by the data. Under genuine prior information, epistemic uncertainty can be represented probabilistically. Under weak or absent prior information, a precise posterior distribution can create a false sense of knowledge. IMs respond by using belief, plausibility, and, in recent work, possibility measures, with validity as the calibration requirement \citep{Martin2019FalseConfidence,Martin2026PossibilisticReview}. 

Because IMs are less familiar than Bayesian posteriors, confidence intervals, or likelihood ratios, this review uses four running examples as a bridge between the abstract construction and practical inference. These examples reappear selectively in later sections so that the reader can see how the IM answer changes when one moves from a basic IM to conditional or marginal IM analyses in a validity-first and efficiency-driven fashion.

\begin{enumerate}[label=\arabic*.]
  \item \textbf{Normal models:} $X_1,\ldots,X_n \stackrel{iid}{\sim} N(\mu,\sigma^2)$. The known-$\sigma$ case gives the cleanest auxiliary pivot and shows how a symmetric PRS reproduces the usual normal confidence interval as an IM plausibility region. It will also be used to separate the calibrated success of the flat Jeffreys prior for $\mu$ from the stronger claim that a formal prior represents objective ignorance. The unknown-$\sigma$ case serves as the simplest example in which conditional IM reasoning and marginal IM reasoning both recover the familiar calibrated inference for $\mu$ without requiring a joint prior for $(\mu,\sigma)$. It also shows why Jeffreys priors can be delicate even in this familiar two-parameter model.
  \item \textbf{Binomial model:} $X\sim\operatorname{Bin}(n,p)$. Following \citet{ClopperPearson1934Binomial}, we consider this example fundamental for developing statistical theory and methods because the inverse association is naturally set-valued, that is, a single auxiliary value can correspond to an interval of $p$ values. That feature makes binomial inference a useful entry point to Dempster-Shafer's set-valued mappings and belief functions, and it shows why random sets and plausibility are more honest than a forced additive posterior-like probability in discrete problems \citep{Dempster1967UpperLower,Shafer1976Evidence}.
  \item \textbf{Eight schools:} the data are the estimated effects of eight SAT coaching programs
  \(
    y=(28,8,-3,7,-1,1,18,12)
  \)
  with known standard errors
  \(
    \sigma=(15,10,16,11,9,11,10,18).
  \)
  In this classic example of Bayesian data analysis \citep{GelmanEtAl2013BDA},
  the usual hierarchical model is
  \begin{equation}\label{eq:eight-school-h}
    y_j\mid\theta_j\sim N(\theta_j,\sigma_j^2),\qquad
    \theta_j\mid\mu,\tau\sim N(\mu,\tau^2),\qquad j=1,\ldots,8.
  \end{equation}
For this eight-schools model, it is better for present purposes first to integrate out the latent school effects in \eqref{eq:eight-school-h} and obtain the marginal sampling model
\begin{equation}\label{eq:eight-school-ms}
   y_j \stackrel{\mathrm{ind}}{\sim} N(\mu,\tau^2+\sigma_j^2), \qquad j=1,\ldots,8.  
   \end{equation} 
  Prior specification for the between-school standard deviation $\tau$ strongly affects the scientific story, especially the amount of pooling and the posterior mass near homogeneity $\tau=0$ \citep{Rubin1981,GelmanEtAl2013BDA,Gelman2006Variance}. In this review the example appears as a constrained problem for $\tau\geq0$ with $\mu$ fixed, illustrating elastic PRSs. Treating $\mu$ as nuisance would require a fuller conditional or marginal IM analysis, which is beyond the scope of this review and will be reported elsewhere. 
  \item \textbf{Bivariate normal correlation:} $(Y_{1i},Y_{2i})$, $i=1,\ldots,n$, are independent bivariate normal observations with known zero means, known unit variances, and unknown correlation coefficient $\rho$. In this example, the sufficient statistic is two-dimensional while the target is one-dimensional. Simpler than the eight-schools example, this is a benchmark for conditional inference because no single global ancillary reduction is available. It motivates local conditional IMs \citep{MartinLiu2015Conditional}.
\end{enumerate}
In addition to these running examples, the Behrens--Fisher problem is introduced later in the discussion of marginal IMs.

The article proceeds historically and conceptually as a review through examples. It starts in Section~\ref{sec:bayes-fiducial} with the background material for IMs: Subsection~\ref{sec:bayes} reviews Bayesian prior specification and questions the probabilistic uncertainty interpretation of information-based priors, while Subsection~\ref{sec:fiducial-confidence} revisits fiducial inference and confidence intervals as prior-free attempts to recover post-data uncertainty. Section~\ref{sec:im} develops the basic IM construction, reviews elastic PRSs for constrained problems, and then works through the running examples. Sections~\ref{sec:conditional} and \ref{sec:marginal} then review two advanced IM developments for efficient inference: conditional IMs reduce auxiliary dimensions by conditioning, while marginal IMs target interest parameters in the presence of nuisance structure. Section~\ref{sec:conclusion} closes with a brief conclusion, points readers to several further unsettled topics, and sheds light on Jaynes' view of probability as the logic of science.

\section{Bayes, Fiducial Inference, and Confidence Intervals}
\label{sec:bayes-fiducial}

\subsection{Bayes and the Problem of Prior Specification}
\label{sec:bayes} 

Bayes requires total probabilistic specification. The difficulty begins when formal priors are presented as objective representations of scientific ignorance. As the most important example, Jeffreys' rule chooses a density proportional to the square root of the Fisher information determinant \citep{Jeffreys1946InvariantPrior,Jeffreys1948TheoryProbability}. 
Logically, Jeffreys' rule solves a particular formal problem: if ignorance could be represented by a density, and if that representation is required to be invariant under smooth changes of coordinates, then the Jeffreys density would be the natural answer. But formal ignorance does not by itself determine an additive probability distribution with calibrated scientific meaning \citep{MartinLiu2015Book}. Jaynes' maximum entropy principle \citep{Jaynes1957InformationTheory,Jaynes2003LogicScience} should be read similarly. Reference priors \citep{Bernardo1979Reference,BergerBernardoSun2009FormalReference} are another principled formal construction, often target- and ordering-dependent, so their usefulness should likewise be separated from the stronger claim that they supply objective prior knowledge.

Although they often produce excellent answers, these priors always require care when interpreted as objective scientific uncertainty. The normal example shows both the success and the danger. If $X\sim N(\mu,\sigma^2)$ with known $\sigma$, the flat prior $\pi(\mu)\propto 1$ gives
\[
  \mu \mid x \sim N(x,\sigma^2)
\]
for one observation, or $N(\bar x,\sigma^2/n)$ for a sample with the sample mean $\bar{x}$. The resulting central credible intervals coincide with the usual confidence intervals and with the fiducial intervals obtained from the pivot $(\bar X-\mu)/(\sigma/\sqrt n)$. While this is a good interval procedure, it does not follow that a flat density is an objective probability distribution over $\mu$, nor that all posterior probabilities computed from it have objective scientific meaning.

With unknown $\sigma$, one must distinguish two different formal priors. The full Jeffreys prior for this two-parameter model, with parameters $(\mu,\sigma)$, is
\[
  \pi_J(\mu,\sigma)\propto \{\det I(\mu,\sigma)\}^{1/2}\propto \frac{1}{\sigma^2}.
\]
The more popular prior is
\[
  \pi(\mu,\sigma)\propto \frac{1}{\sigma}
\]
which is the right-Haar, independence Jeffreys, and reference prior for the usual location-scale analysis. It is this modified Jeffreys prior that matches the standard $t$ posterior for $\mu$ after integrating out $\sigma$. Again, the resulting one- and two-sided intervals have a compelling confidence or plausibility interpretation. But that success belongs to the pivot and its calibration, not to the claim that the full Jeffreys prior $\sigma^{-2}$, or any other formal prior, is an objective probability representation of ignorance over $(\mu,\sigma)$. 

For $p\in (0,1)$ in the binomial example $X\sim\operatorname{Bin}(n,p)$, Jeffreys' prior is $\operatorname{Beta}(1/2,1/2)$, while a uniform formal prior \citep{Bayes1763Essay} gives $\operatorname{Beta}(1,1)$. These choices lead to different posterior intervals, especially near $x=0$ or $x=n$.
The bivariate normal correlation coefficient gives a different warning about whether formal priors can be consistent with Bayes as a rule for sequentially combining information. The parameter $\rho\in(-1,1)$ is scalar, but the natural sufficient statistic has two components. To see the ambiguity, transform the observations to
\[
  Y'_{1i}=2^{-1/2}(Y_{1i}+Y_{2i}),\qquad
  Y'_{2i}=2^{-1/2}(Y_{1i}-Y_{2i}),
\]
so that the two transformed samples are independent with
\[
  Y'_{1i}\sim N(0,1+\rho),\qquad
  Y'_{2i}\sim N(0,1-\rho).
\]
If the $Y'_{1i}$'s are regarded as the first data set, Jeffreys' rule for this first marginal model gives
\[
  \pi_+(\rho)\propto \frac{1}{1+\rho}.
\]
A sequential Bayesian analysis then updates this prior by the likelihood for the $Y'_{1i}$'s and uses the resulting posterior as the prior when the $Y'_{2i}$'s become available; the final posterior is therefore proportional to the joint likelihood times $\pi_+(\rho)$. Reversing the order gives instead the first-stage Jeffreys prior
\[
  \pi_-(\rho)\propto \frac{1}{1-\rho},
\]
and the final posterior is proportional to the same joint likelihood times $\pi_-(\rho)$. A one-shot Jeffreys calculation from the joint transformed model gives a third prior,
\[
  \pi_J(\rho)\propto \frac{(1+\rho^2)^{1/2}}{1-\rho^2}.
\]
Thus the sequential-updating route depends on which transformed component is processed first, while the joint-model route gives a different formal prior again. 

The eight-schools example is even more revealing. The parameter $\tau$ controls heterogeneity. A prior on $\tau$ is not a harmless technicality. It determines how strongly the analysis shrinks school effects toward a common mean, and how much posterior mass is placed near homogeneity. The formal-prior ambiguity is similar to both the normal example and the bivariate normal correlation example: applying Jeffreys or reference-prior reasoning to the marginal model for $(\mu,\tau)$. Flat priors, half-Cauchy priors, and Jeffreys-type priors can therefore give noticeably different scientific narratives. This is not a failure of Bayes when the prior is genuine; it is a warning against calling formal priors objective.

The position taken here is therefore sharp but not anti-Bayesian. If one wants to do Bayes, one should be genuinely Bayesian and specify the prior as real scientific, subjective, or empirical. If one does not have such information, then IMs offer a way to quantify uncertainty without pretending that a formal prior has solved the problem. Information-based Bayes remains indispensable, especially for understanding probability as a disciplined logic of uncertain reasoning. But the scientific interpretation of its posterior probabilities must be earned, often through calibration, matching, or an IM plausibility interpretation, rather than assumed from the formal prior alone. This is close in spirit to Rubin's applied Bayesian position that frequency calculations can be essential for understanding, communicating, validating, and monitoring Bayesian analyses \citep{Rubin1984Bayesianly}.

\subsection{Fiducial Inference and Confidence Intervals}
\label{sec:fiducial-confidence}

Historically, Fisher's fiducial argument preceded Neyman's mature confidence theory. Expositionally, however, it is helpful to begin with familiar confidence intervals, because Neyman's construction gives the mathematically clean interpretation of many fiducial intervals. Neyman's 1941 paper is explicitly titled ``Fiducial Argument and the Theory of Confidence Intervals,'' and the later confidence-distribution literature also treats many fiducial distributions as confidence distributions in regular one-parameter settings \citep{Neyman1941FiducialConfidence,XieSingh2013ConfidenceDistribution,SchwederHjort2016}. 

In Neyman's theory, probability attaches to the random procedure before the data are observed \citep{Neyman1937Confidence}. A $100(1-\alpha)\%$ confidence procedure $C(X)$ satisfies
\[
  \pr_\theta\{\theta\in C(X)\} \geq 1-\alpha
\]
for each $\theta$ in the parameter space. Once $X=x$ is observed, $C(x)$ is a fixed set. The classical theory does not say that the fixed $\theta$ lies in this fixed set with probability $1-\alpha$. The gain is conceptual clarity; the cost is a limited post-data evidential interpretation.

Fisher wanted more for situation-specific probabilistic inference. In the known-variance normal model,
\[
  \bar X = \mu + \frac{\sigma}{\sqrt n}Z, \qquad Z\sim N(0,1),
\]
the observed value $\bar x$ of the sample mean $\bar{X}$ gives the formal inversion
\[
  \mu = \bar x - \frac{\sigma}{\sqrt n} Z.
\]
If one continues to regard $Z$ as standard normal after $\bar X=\bar x$ is observed, then $\mu$ receives the fiducial distribution $N(\bar x,\sigma^2/n)$. Its central intervals match the usual confidence intervals. In this example, Fisher's fiducial answer, Neyman's confidence interval, and the flat-prior Bayesian interval agree numerically. The unknown-variance normal model gives the next classical case. The pivot
\begin{equation}\label{eq:normal-example-T}
   T=\frac{\bar X-\mu}{S/\sqrt n}\sim t_{n-1},  
\end{equation}
where $S^2$ is the unbiased sample variance as the estimator of $\sigma^2$,
leads to the usual $t$ interval for $\mu$. It can be read as a confidence interval, a fiducial interval, or as a formal Bayesian credible interval under the right-Haar/reference prior $\pi(\mu,\sigma)\propto 1/\sigma$.

The binomial case shows why the above equivalence breaks down. If $X\sim\operatorname{Bin}(n,p)$, a natural auxiliary representation uses $U\sim\Unif(0,1)$ and
\[
  X = F_p^{-1}(U),
\]
where $F_p$ is the binomial distribution function. After observing $X=x$, the inversion is not a single value of $p$ for each $u$; it is generally a set determined by
\begin{equation}\label{eq:binomial-example-set}
  F_p(x-1) \leq u < F_p(x).
\end{equation}
\citet{ClopperPearson1934Binomial} explicitly presented their construction in confidence or fiducial language for this binomial problem. Its intervals are exact but conservative, and discreteness exposes the gap between calibrated set inference and precise posterior-like probability.

The bivariate normal correlation coefficient and eight schools expose a second gap, this time caused by conditioning. In the bivariate normal correlation, for example, the transformed sums
\[
  S_1=\frac12\sum_{i=1}^n(Y_{1i}+Y_{2i})^2,\qquad
  S_2=\frac12\sum_{i=1}^n(Y_{1i}-Y_{2i})^2
\]
serve as sufficient statistics and satisfy $S_1=(1+\rho)U_1$ and $S_2=(1-\rho)U_2$, where $U_1,U_2$ are independent $\chi^2_n$ variables. It is tempting to make fiducial inference by conditioning on an identity such as $S_1/U_1+S_2/U_2=2$. But that identity is not a separable ancillary reduction of the form needed for valid conditioning. Classical likelihood theory handles this through ancillary or approximate conditional arguments, including the $r^\star$ machinery of \citet{BarndorffNielsen1986SignedLogLikelihood}, \citet{Fraser1990ObservedLikelihood}, and \citet{Reid1995Conditioning}. The IM response is to ask a more primitive question: {\it which lower-dimensional auxiliary quantity can be predicted validly after the observed auxiliary features have been conditioned on?} Conditional IMs answer this question later in Section \ref{sec:conditional}.

Nevertheless, fiducial reasoning has remained a source of inspiration and
the historical study of \citet{Zabell1992FisherFiducial} serves as an essential guide. Fisher identified a real scientific need for posterior-like reasoning without arbitrary priors. But the fiducial argument did not become a general coherent probability calculus. 
Fraser's structural equations \citep{Fraser1968Structure}, Dempster-Shafer theory \citep{Dempster1967UpperLower,Shafer1976Evidence}, and generalized fiducial inference \citep{Hannig2009GFI,HannigIyerLaiLee2016GFIReview} are three later developments inspired by this unfinished ambition. These are crucial predecessors of IMs. The IM contribution is to keep the auxiliary-variable source of uncertainty but change the order of operations and seek ways to improve efficiency through conditioning and marginalization in the auxiliary space.

\section{Inferential Models}
\label{sec:im}

\subsection{The Basic IM Construction}
\label{sec:im-basic}

The IM framework begins with an association
\begin{equation}\label{eq:association}
    X = a(\theta,U), \qquad U\sim\pr_U,
\end{equation}
where $X$ is observable, $\theta$ is unknown, and $U$ is an auxiliary variable with known distribution. After $X=x$ is observed, there exists an unobserved auxiliary value $u^\star$ consistent with the true parameter and the observed data. The inferential problem is to reason about $\theta$ through uncertainty about $u^\star$.

The distinction from fiducial reasoning is the central point. Fiducial methods transfer the auxiliary variable to the parameter space first, producing a data-dependent parameter distribution, and then apply probability calculus in the parameter space. IMs carry out the probability operation before that transfer \citep{LiuMartin2015WIREs}; see also \citet{Barnard1995Pivotal}. They predict $u^\star$ in the well-defined auxiliary-variable probability space by a valid PRS, and only after this prediction do they combine the result with the observed association.

Procedurally, the basic IM has three steps \citep{MartinLiu2013IM,MartinLiu2015Book}.

\begin{enumerate}[label=\arabic*.]
  \item \textbf{Association.} Specify the association \eqref{eq:association} that preserves the sampling model of $X$. For fixed $x$ and $u$, this defines the compatible parameter set
  \[
    \Theta_x(u)=\{\theta:x=a(\theta,u)\}.
  \]
  \item \textbf{Prediction.} Predict the unobserved $u^\star$ by a PRS $\mathcal S\sim\pr_{\mathcal S}$ on the auxiliary space $\mathbb U$. The mathematical calibration requirement for $\mathcal S$ is imposed before any transfer to the parameter space.
  \item \textbf{Combination.} Map the predictive random set to the parameter space:
  \[
    \Theta_x(\mathcal S)=\bigcup_{u\in\mathcal S}\Theta_x(u).
  \]
\end{enumerate}

For an assertion $A\subseteq\Theta$, define the belief and plausibility
\[
  \bel_x(A)=\pr_{\mathcal S}\{\Theta_x(\mathcal S)\subseteq A\},
  \qquad
  \pl_x(A)=\pr_{\mathcal S}\{\Theta_x(\mathcal S)\cap A\neq\emptyset\}.
\]
Thus, belief measures definite support for $A$; plausibility measures the extent to which $A$ has not been ruled out. 

This basic IM is oriented toward validity first. Formally, the calibration logic can be stated as two definitions and a theorem.

\begin{definition}[Valid PRS]\label{def:valid-prs}
For a PRS $\mathcal S\sim\pr_{\mathcal S}$ on the auxiliary space $\mathbb U$, define the PRS contour
\begin{equation}\label{eq:prs-contour}
  \gamma_{\mathcal S}(u)=\pr_{\mathcal S}\{\mathcal S\ni u\},\qquad u\in\mathbb U.
\end{equation}
The PRS is \emph{valid} for predicting the auxiliary value if
\begin{equation}\label{eq:prs-validity}
  \pr_U\{\gamma_{\mathcal S}(U)\leq \alpha\}\leq \alpha,\qquad 0<\alpha<1,
\end{equation}
where $U\sim\pr_U$. Equivalently, $\gamma_{\mathcal S}(U)$ is stochastically no smaller than a $\Unif(0,1)$ random variable. In miss-probability form, with
\[
  Q_{\mathcal S}(u)=\pr_{\mathcal S}\{\mathcal S\not\ni u\}=1-\gamma_{\mathcal S}(u),
\]
validity is
\[
  \pr_U\{Q_{\mathcal S}(U)\geq 1-\alpha\}\leq \alpha,\qquad 0<\alpha<1.
\]
Thus a large miss probability can occur only on a small $\pr_U$-set of auxiliary values.
\end{definition}

\begin{definition}[Valid posterior assertion plausibility]\label{def:valid-plausibility}
For an assertion $A\subseteq\Theta$, the IM procedure gives the posterior assertion plausibility
\[
  \pl_x(A)=\pr_{\mathcal S}\{\Theta_x(\mathcal S)\cap A\neq\emptyset\}.
\]
This plausibility is \emph{valid} for assertion $A$ if
\begin{equation}\label{eq:pl-validity}
  \sup_{\theta\in A}\pr_{X|\theta}\{\pl_X(A)\leq \alpha\}\leq \alpha,\qquad 0<\alpha<1.
\end{equation}
Equivalently, when $A$ is true, $\pl_X(A)$ is stochastically no smaller than a $\Unif(0,1)$ random variable. The IM is valid if \eqref{eq:pl-validity} holds for all assertions $A\subseteq\Theta$.
\end{definition}

The basic IM validity theorem of \citet{MartinLiu2013IM} links auxiliary prediction to assertion-wise calibration.

\begin{theorem}[Validity theorem]\label{thm:im-validity}
If the PRS $\mathcal S$ is valid in the sense of Definition~\ref{def:valid-prs} and $\Theta_x(\mathcal S)$ is non-empty with $\pr_{\mathcal S}$-probability one for each $x$, then the IM obtained by the association, prediction, and combination steps is valid in the sense of Definition~\ref{def:valid-plausibility} for every assertion $A$. Consequently, the plausibility region
\[
  \Pi_x(\alpha)=\{\theta:\pl_x(\{\theta\})>\alpha\}
\]
satisfies
\[
  \inf_{\theta\in\Theta}\pr_{X|\theta}\{\theta\in \Pi_X(\alpha)\}\geq 1-\alpha,
\]
so the fixed-data plausibility scale has a direct frequentist calibration.
\end{theorem}

\subsection{Elastic PRSs for Constrained Problems}
\label{sec:elastic-prs}

The basic construction above assumes that $\Theta_x(\mathcal S)$ is non-empty with probability one. This assumption can fail, for example, in constrained problems. Suppose the parameter is known to lie in a proper subset $C\subseteq\Theta$, such as a nonnegative mean, a variance component, or a Poisson signal rate with known background. For fixed data $x$, the constraint maps back to the auxiliary space through
\[
  U_{C,x}=\bigcup_{\theta\in C}\{u:x=a(\theta,u)\}.
\]
If the realized PRS $\mathcal S$ misses $U_{C,x}$, then the corresponding parameter random set is incompatible with the constraint. These are referred to as conflict cases.

One solution to the conflict case problem is to condition away the conflict cases, as in Dempster's rule of combination. This is valid in many settings, but it can be inefficient because it discards informative conflict cases or redistributes more predictive mass than is needed. Elastic PRSs were introduced by \citet{ErminiLeafLiu2012ElasticBelief} to remove conflict cases more locally. Start with a valid baseline PRS $\mathcal S_0$ and embed it in an increasing family
\[
  \mathcal S_e,\qquad 0\leq e\leq 1,
\]
where $\mathcal S_0$ is the original PRS and $\mathcal S_1$ is large enough to guarantee compatibility. For each realized auxiliary prediction, choose the smallest stretching amount
\[
  \widehat e=\inf\{e:\mathcal S_e\cap U_{C,x}\neq\emptyset\},
\]
and combine with
\[
  \Theta_x(\mathcal S_{\widehat e})\cap C.
\]
Thus the PRS stretches only until it intersects the auxiliary values compatible with the constraint. The elastic step is data-dependent, but it is not an arbitrary posterior truncation. It is carried out in the auxiliary space, where validity is controlled, and then transferred to the constrained parameter space.

The key property is that this minimal stretching preserves calibration while improving efficiency. \citet{ErminiLeafLiu2012ElasticBelief} show that the resulting elastic belief function is valid for assertions $A\subseteq C$ and is typically more efficient than Dempster's rule of combination.  Elastic PRSs are important conceptually because constraints are common in scientific inference. Variance components, signal rates, scale parameters, probabilities, and physical effects often live on restricted spaces. Bayesian methods usually encode such restrictions by truncating or reshaping a prior.

\subsection{Illustrative Examples}
\label{sec:im-examples}

\begin{figure}[tbp]
  \centering
  \includegraphics[width=\textwidth]{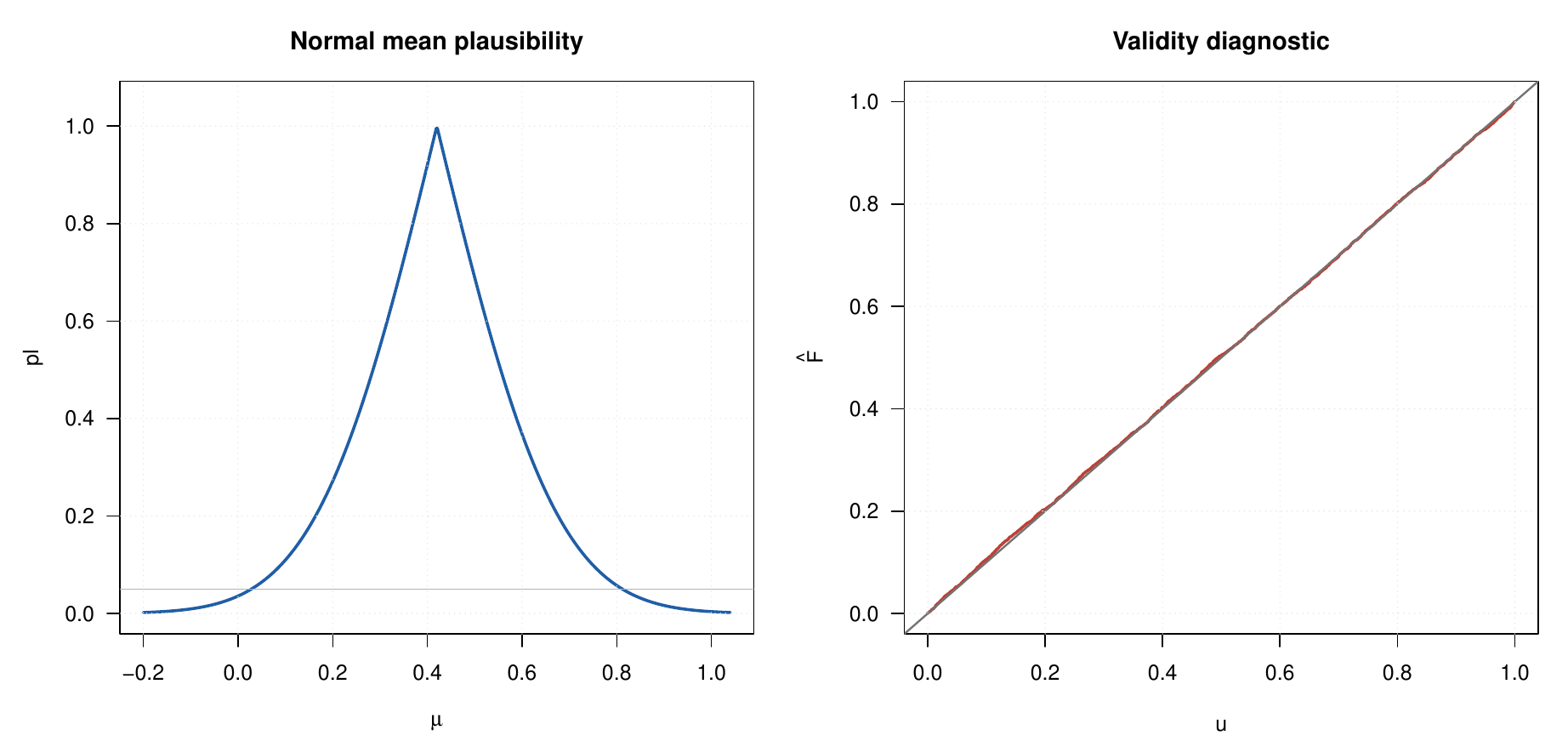}
  \caption{Normal mean example. The left panel shows the plausibility curve for $\mu$ when $n=25$, $\sigma=1$, and $\bar x=0.42$. The right panel is a validity diagnostic based on 5000 simulations under $\mu=0$; the empirical distribution of $\pl_X(0)$ closely follows the uniform diagonal.}
  \label{fig:normal-plausibility-validity}
\end{figure}

\begin{figure}[tbp]
  \centering
  \includegraphics[width=\textwidth]{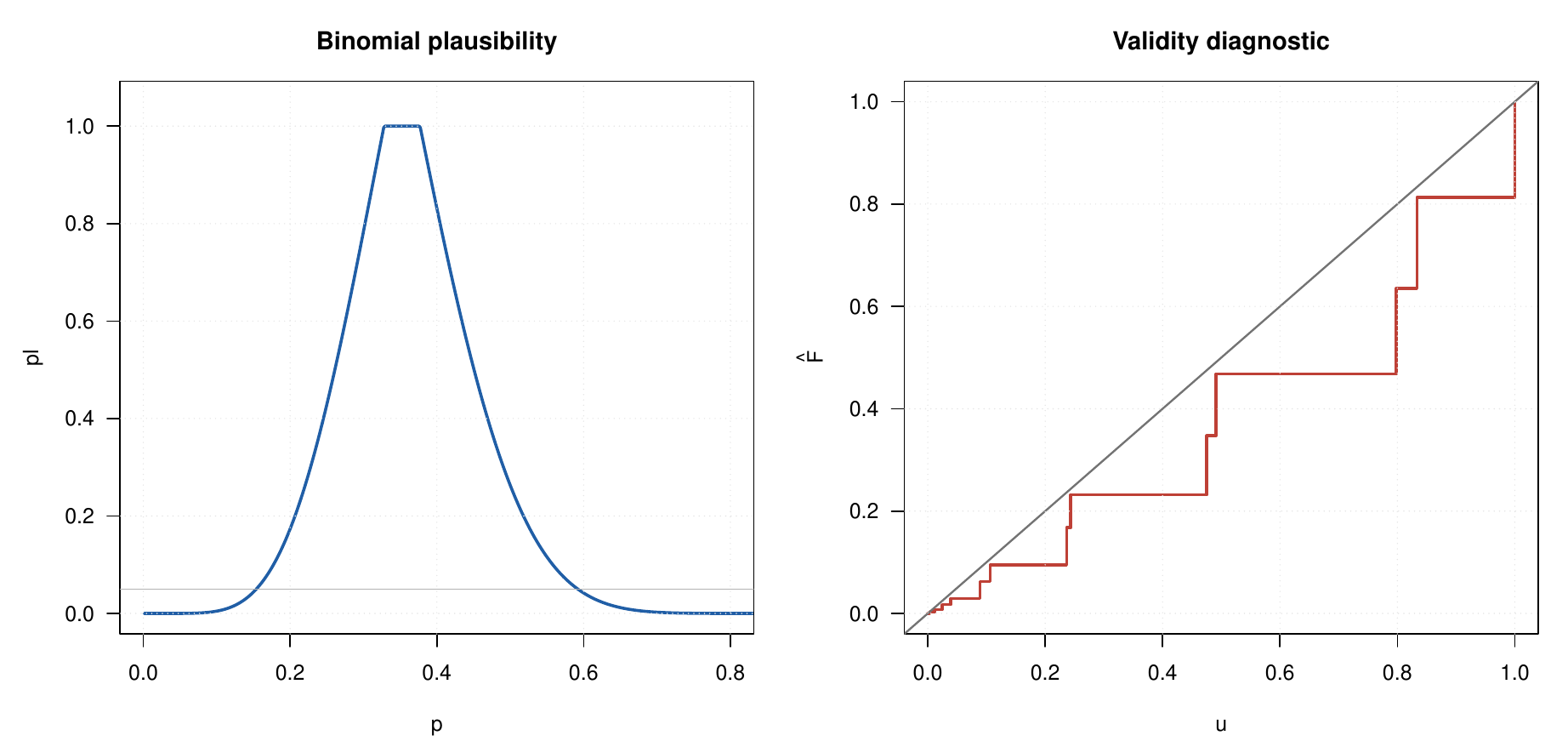}
  \caption{Binomial example. The left panel shows an exact two-sided plausibility curve for $p$ when $n=20$ and $x=7$. The right panel shows a validity diagnostic based on 5000 simulations under $p=0.35$; discreteness makes the distribution conservative, so the empirical curve lies mostly below the uniform diagonal.}
  \label{fig:binomial-plausibility-validity}
\end{figure}

\begin{figure}[tbp]
  \centering
  \includegraphics[width=\textwidth]{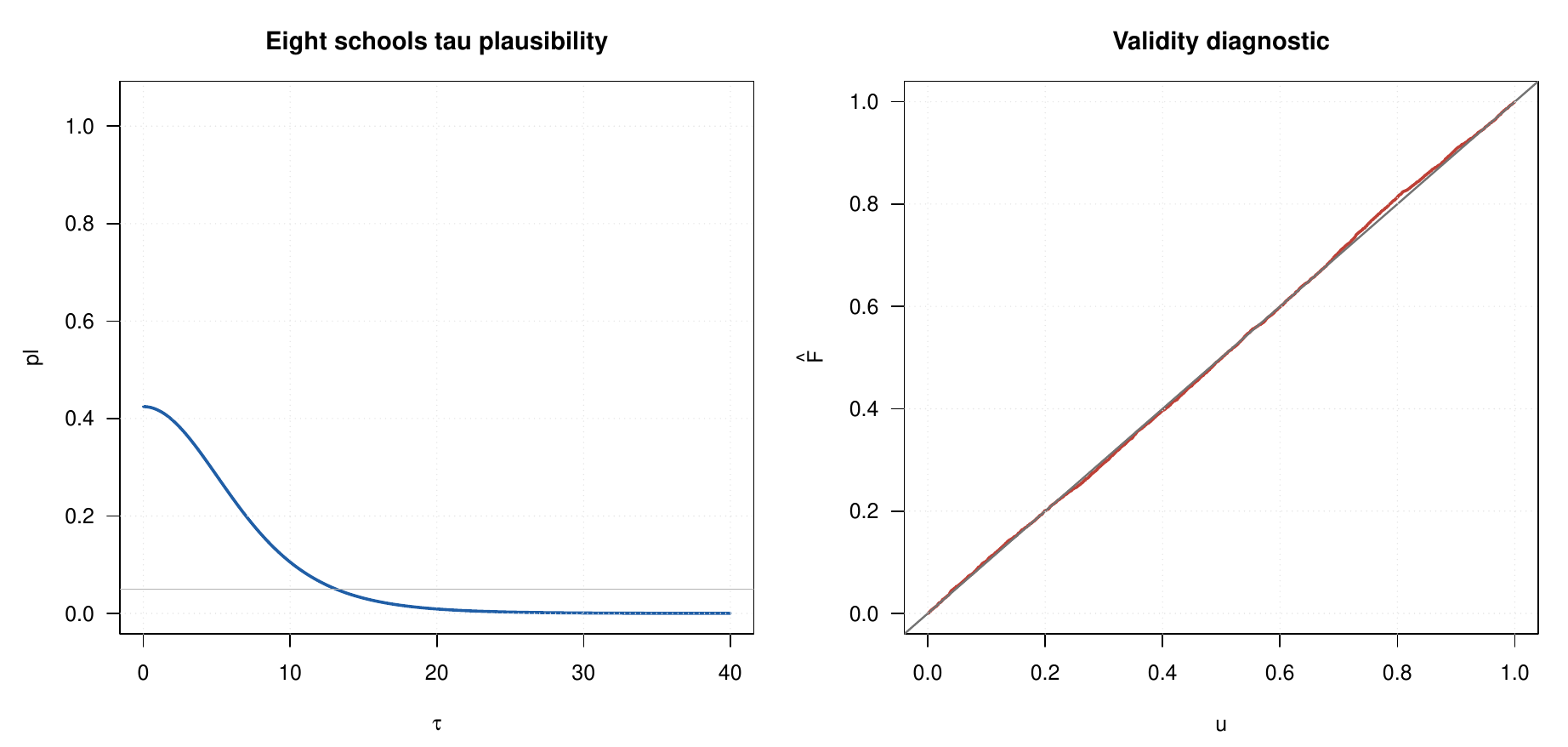}
  \caption{Eight-schools heterogeneity example with $\mu=8$ fixed. The left panel shows the plausibility curve for the constrained parameter $\tau\geq0$ using the marginal model after integrating out the school effects. The curve is largest at the boundary, reflecting little evidence for extra heterogeneity under this simple known-$\mu$ illustration. The estimated probability that the baseline central PRS is in conflict and therefore needs elastic stretching is $0.576$. The right panel is a validity diagnostic based on 5000 simulations under $\tau=10$.}
  \label{fig:eight-schools-plausibility-validity}
\end{figure}

The following selected numerical illustrations pair fixed-data plausibility curves with simulation diagnostics for validity. In the diagnostics, the plausibility assigned to the true parameter value is computed repeatedly under the model and compared with the $\Unif(0,1)$ benchmark; exact continuous IMs should track the diagonal, while discrete examples often unavoidably appear conservative. The diagnostic plots are shown in Figures 
\ref{fig:normal-plausibility-validity},
\ref{fig:binomial-plausibility-validity},
and
\ref{fig:eight-schools-plausibility-validity}
for the normal, binomial, and eight-school examples that are discussed below.

In the known-variance normal example, the conditional-IM reduction (see Section~\ref{sec:conditional}) gives the association for the reduced sufficient statistic
\[
  \bar X=\mu+\frac{\sigma}{\sqrt n}Z, \qquad Z\sim N(0,1).
\]
An estimation-oriented default symmetric predictive random set is
\[
  \mathcal S=\{z: |z|\leq |Z'|\}, \qquad Z'\sim N(0,1).
\]
For the singleton assertion $\{\mu\}$, the plausibility contour is
\[
  \pl_x(\mu)=1-\left|2\Phi\left(\frac{\bar x-\mu}{\sigma/\sqrt n}\right)-1\right|.
\]
The plausibility set $\{\mu:\pl_x(\mu)>\alpha\}$ is exactly the usual two-sided $100(1-\alpha)\%$ normal interval. Thus the familiar interval is not discarded. It is reinterpreted as a calibrated plausibility region arising from valid prediction of the auxiliary variable and giving the desired interpretation of situation-specific inference.

For the unknown-variance normal problem, after efficiency operations (see Sections \ref{sec:conditional} and \ref{sec:marginal}),
one can work with \eqref{eq:normal-example-T}.
Using the same symmetric predictive-random-set idea on the $t$ auxiliary scale gives
\[
  \pl_x(\mu)=1-\left|2F_{t_{n-1}}\left(\frac{\bar x-\mu}{s/\sqrt n}\right)-1\right|,
\]
and the plausibility regions coincide with the usual $t$ intervals. Again, this agreement does not say that all formal posterior probabilities are objective. It says that the interval has a calibrated auxiliary-variable predictive inference justification.

For the binomial model, the association $X=F_p^{-1}(U)$ with $U\sim\Unif(0,1)$ gives \eqref{eq:binomial-example-set} and
\[
  \Theta_x(u)=\{p:F_p(x-1)\leq u<F_p(x)\}.
\]
Because $\Theta_x(u)$ is set-valued and the distribution is discrete, a singleton fiducial transfer is unnatural. Like Dempster-Shafer, an IM keeps the set-valued nature visible. Predicting $u^\star$ by a valid PRS yields plausibility regions closely related to exact binomial tests and Clopper-Pearson intervals, but the output is belief and plausibility rather than a forced confidence interval for $p$.

For the eight-schools model illustration, take $\mu$ to be known and focus on $\tau\geq 0$. This turns the example into a clean constrained-parameter problem. An auxiliary association is
\[
  Y_j=\mu+(\tau^2+\sigma_j^2)^{1/2}Z_j,\qquad Z=(Z_1,\ldots,Z_8)\sim N_8(0,I).
\]
For a candidate $\tau$, the auxiliary value required by the observation is
\[
  z_\tau(y)=\left(\frac{y_1-\mu}{(\tau^2+\sigma_1^2)^{1/2}},\ldots,
  \frac{y_8-\mu}{(\tau^2+\sigma_8^2)^{1/2}}\right).
\]
Given a PRS $\mathcal S$ for $Z$, the value $\tau$ is plausible if $z_\tau(y)\in\mathcal S$. Without further reducing dimension via conditional IMs, for example, a central nested PRS based on $h(z)=\sum_j z_j^2$ gives a plausibility contour of the form
\[
  \pl_y(\tau)=1-\left|2F_{\chi^2_8}\{h(z_\tau(y))\}-1\right|,
  \qquad \tau\geq 0,
\]
where $F_{\chi^2_8}$ is the $\chi^2_8$ distribution function. This two-sided form treats both unusually large and unusually small standardized residual sums as atypical. The constraint $\tau\geq 0$ is not cosmetic: if the combined random set $\{\tau\geq0:z_\tau(y)\in\mathcal S\}$ is empty for a realized PRS, the elastic PRS construction stretches the auxiliary prediction only enough to hit the curve $\{z_\tau(y):\tau\geq0\}$. For the observed data with $\mu=8$, $h(z_0(y))=4.713$, $\pl_y(0)=0.424$, and the probability that the baseline central PRS misses the constrained auxiliary curve is $1-\sup_{\tau\geq0}\pl_y(\tau)=0.576$. This is the probability that the elastic adjustment is activated in this illustration.

This conflict probability is not itself a point on the plausibility curve. The curve reports, for each fixed $\tau$, how plausible that value is after combining the data, the PRS, and the constraint. The conflict probability is a diagnostic of the baseline PRS before elastic repair: it measures how often the auxiliary prediction fails to hit any point on the constrained curve $\{z_\tau(y):\tau\geq0\}$. Thus it is naturally reported as a separate additional annotation or diagnostic. The case with unknown $\mu$ introduces a nuisance parameter and is better treated as both a conditional IM and a marginal IM problem. A full treatment thus requires substantial space and will be reported elsewhere.

The bivariate normal correlation example is not developed here because its main role is to motivate local conditional IMs, where it will be treated explicitly. The philosophical lesson across the examples is the same. The sampling model supplies an auxiliary source of randomness. Unlike fiducial reasoning, IMs keep the randomness in its native auxiliary space, calibrate the prediction there, and then transfer the resulting set-valued uncertainty to the parameter space. This is why IMs can preserve Fisher's post-data ambition and Neyman's calibration requirement at the same time.

\section{Conditional IMs}
\label{sec:conditional}

Conditional IMs  \citep{MartinLiu2015Conditional} are the first major efficiency refinement of the basic IM construction. A baseline association may contain many auxiliary coordinates even when the parameter or assertion of interest is low-dimensional. Predicting all of those coordinates is valid but wasteful, and wasteful prediction produces inefficient plausibility. The conditional IM idea is to identify auxiliary features that have effectively been observed, condition on them, and predict only the remaining unobserved auxiliary feature. Classical reviews of ancillary conditioning make clear why this matters: conditioning can define the relevant subset of the sample space and can make inference more efficient \citep{Fraser2004Ancillaries,Reid1995Conditioning,GhoshReidFraser2010Ancillary}. The difference is that conditional IMs carry out conditioning in the auxiliary-variable space.

\subsection{Regular Conditional IMs}
\label{sec:regular-cims}

Start from a baseline association \eqref{eq:association}.
A regular conditional IM is available when the association can be rewritten, after one-to-one transformations of data and auxiliary variables, as
\begin{equation}\label{eq:cim-association}
   h_1(X)=a_1\{\phi_1(U),\theta\},\qquad
  h_2(X)=a_2\{\phi_2(U)\},   
\end{equation}
where the second equation is free of $\theta$. The observed value $h=h_2(x)$ reveals the corresponding feature of the auxiliary variable. Therefore, instead of predicting the full $U$, the conditional IM predicts
\[
  V=\phi_1(U)
\]
under the conditional distribution of $V$ given $a_2\{\phi_2(U)\}=h$. The association used for inference is the lower-dimensional conditional association
\[
  h_1(x)=a_1(v,\theta),\qquad V\sim\pr_{\phi_1(U)\mid a_2\{\phi_2(U)\}=h}.
\]
The A/P/C steps that produce inferential assessments are then the same as before: invert this conditional association, predict the unobserved $V^\star$ by a predictive random set on the conditional auxiliary space, and combine.

The important change is that validity is now conditional. If the predictive random set $\mathcal S$ is conditionally valid under the conditional law at $h$, then the resulting plausibility is calibrated on the relevant subset $\{x:h_2(x)=h\}$. In particular, for a true assertion $A$,
\[
  \pr_\theta\{\pl_X(A)\leq\alpha\mid h_2(X)=h\}\leq \alpha,
  \qquad 0<\alpha<1,
\]
for $\theta\in A$
under the usual non-empty-focal-set condition. This is stronger and more interpretable than an unconditional guarantee when the conditioning statistic really is part of the observed auxiliary information.

Regular conditional IMs clarify the relation between IMs, sufficiency, ancillarity, and Bayes. If $h_1(X)$ is sufficient and the conditional law of $\phi_1(U)$ does not depend on $h$, then the conditional reduction reproduces the familiar instruction to work with the sufficient statistic. If $h_2(X)$ is ancillary, then validity is assessed conditionally on the ancillary value, matching Fisher's relevant-subset idea. From the IM viewpoint, both are special cases of a more primitive operation: combine information by conditioning on observed auxiliary features before prediction. This is also why conditional IMs can be viewed as a prior-free counterpart to Bayesian information aggregation, but without requiring a prior distribution on $\theta$.

The word ``regular'' is doing real work. The same conditioning map must apply across the parameter space, and the conditioning equation \eqref{eq:cim-association} must separate data from auxiliary variables in a parameter-free way. Many textbook models have this form after a sensible transformation. But there are important models in which the transformation \eqref{eq:cim-association} does not exist. In such cases, no single global conditioning feature is available from \eqref{eq:association}, and one can consider alternative specifications of \eqref{eq:association} or local conditional IMs.

\subsection{Local Conditional IMs}
\label{sec:local-cims}

Local conditional IMs \citep{MartinLiu2015Conditional} were introduced for precisely those nonregular cases. The idea is to assess each singleton assertion $\{\theta_0\}$ using a conditional association localized at $\theta_0$. For a fixed candidate value, one seeks transformations
\[
  T_{\theta_0}(X)=b_{\theta_0}\{\tau_{\theta_0}(U),\theta_0\},\qquad
  H_{\theta_0}(X)=\eta_{\theta_0}(U),
\]
at least locally around $\theta_0$, where $H_{\theta_0}(X)$ is the observed auxiliary feature to condition on and $\tau_{\theta_0}(U)$ is the reduced auxiliary variable to be predicted. The plausibility of $\theta_0$ is computed from the conditional distribution
\[
  \tau_{\theta_0}(U)\mid \eta_{\theta_0}(U)=H_{\theta_0}(x).
\]
Then the local plausibility region is assembled as
\[
  \{\theta:\pl_{x,\theta}(\{\theta\})>\alpha\},
\]
where the subscript reminds us that the conditional association used to test $\theta$ or compute $\mbox{pl}_x(\{\theta\})$ is localized at that same value.

This may look unusual at first because the conditioning law changes with the candidate value. But that is exactly what prevents an invalid fiducial shortcut. Local conditional IMs pay a small price in notation and computation to preserve the essential IM ordering: first identify the auxiliary feature still needing prediction, then predict it validly, and only then transfer uncertainty to the parameter assertion.

The validity statement is local but still meaningful. If, at $\theta_0$, the conditional PRS is valid under the conditional law used by the localized association, then
\[
  \pr_{\theta_0}\{\pl_{X,\theta_0}(\{\theta_0\})\leq\alpha
    \mid H_{\theta_0}(X)=h\}\leq\alpha.
\]
Thus a local conditional IM supplies a calibrated plausibility scale for each candidate value, even when no global ancillary or sufficient reduction is available. 

\subsection{A differential geometric formulation}
\label{sec:diffgeo-cims}

The construction can be described geometrically \citep{Liu2026DifferentialGeometricIMs}. For a smooth association $x=a(\theta,u)$, fixing $x$ and varying $\theta$ traces a curve or surface through the auxiliary space. Assume, locally, that the derivative (square matrix) $D_u a(\theta,u)$ is nonsingular. Differentiating the identity $x=a(\theta,u_{x,\theta})$ with $x$ fixed gives the parameter-motion vector fields
\[
  G_j^\theta(u) \equiv \frac{\partial u_{x, \theta}}{\partial \theta_j}
  =-\{D_u a(\theta,u)\}^{-1}\partial_{\theta_j}a(\theta,u),
  \qquad j=1,\ldots,p.
\]
At a candidate value $\theta_0$, a local conditioning map $\eta_{\theta_0}$ should be insensitive to these motions. In differential form it solves
\begin{equation}\label{eq:local-cim-pde}
  D\eta_{\theta_0}(u)G_j^{\theta_0}(u)=0,
  \qquad j=1,\ldots,p.
\end{equation}
Thus local conditioning variables are first integrals of the so-called {\it distribution} in differential geometry
\[
  \mathcal D_{\theta_0}(u)
  =\operatorname{span}\{G_1^{\theta_0}(u),\ldots,G_p^{\theta_0}(u)\}
\]
on the auxiliary space. 

This is the differential formulation behind local conditional IMs. It includes the class of regular conditional IMs as a special case.
A precise mathematical theory is developed in \citep{Liu2026DifferentialGeometricIMs}.

\subsection{Illustrative Examples}
\label{sec:conditional-examples}

\subsubsection{Normal location and scale}
The normal location model shows regular conditioning in its simplest form. With $X_i=\theta+U_i$ and $U_i\stackrel{\mathrm{iid}}{\sim}N(0,1)$,
\[
  \bar X=\theta+\bar U,\qquad X_i-\bar X=U_i-\bar U,\quad i=1,\ldots,n.
\]
The residual vector is an observed auxiliary feature, and $\bar U$ is independent of the residuals. The conditional IM therefore predicts only $\bar U\sim N(0,1/n)$, giving the same plausibility contour as the usual normal pivot. This is a regular conditional IM because the conditioning feature does not depend on $\theta$.

The unknown-variance normal model illustrates how conditional and marginal thinking meet. The details are discussed in Section \ref{sec:marginal}.

\subsubsection{Binomial model}

This example connects directly to Dempster's 1960s innovation of lower and upper probabilities. Instead of starting from the scalar inverse-cdf representation, take
\[
  U_1,\ldots,U_n\stackrel{\mathrm{iid}}{\sim}\Unif(0,1),
  \qquad X_i=1_{\{U_i\leq p\}},\qquad X=\sum_{i=1}^n X_i .
\]
After observing only $X=x$, the individual labels of the uniforms carry no information about $p$. If $U_{(1)}\leq\cdots\leq U_{(n)}$ are the ordered uniforms, with the conventions $U_{(0)}=0$ and $U_{(n+1)}=1$, then the auxiliary values compatible with the count determine the random interval
\[
  \Gamma_x(U)=\{p:U_{(x)}\leq p\leq U_{(x+1)}\},
\]
up to endpoint conventions that are immaterial under the continuous auxiliary law. Dempster's multivalued mapping idea assigns lower and upper probabilities by
\[
  \bel_x(A)=P\{\Gamma_x(U)\subseteq A\},\qquad
  \pl_x(A)=P\{\Gamma_x(U)\cap A\neq\emptyset\},
\]
for assertions $A\subseteq[0,1]$ \citep{Dempster1966NewMethods,Dempster1967UpperLower,Dempster1968RandomClosedInterval}. For example,
\[
  \pl_x\{p\leq p_0\}=P\{U_{(x)}\leq p_0\}=P_{p_0}(X\geq x) \quad\mbox{and}\quad
  \pl_x\{p\geq p_0\}=P\{U_{(x+1)}\geq p_0\}=P_{p_0}(X\leq x).
\]
Thus the ordered-uniform construction recovers the exact binomial tail probabilities as plausibilities, while retaining the set-valued inverse image forced by discreteness.

In contrast, the IM solution can be made even more direct by using the probability integral transform for the observed count itself. Let $U\sim\Unif(0,1)$ and generate
\[
  X=F_p^{-1}(U),
\]
where $F_p$ is the binomial cdf. After observing $X=x$, a candidate $p$ is compatible with an auxiliary value $u$ exactly when
\[
  F_p(x-1)<u\leq F_p(x),
\]
with the usual convention $F_p(-1)=0$. Thus inference for $p$ can be based on a one-dimensional auxiliary variable rather than on Dempster's ordered-uniform endpoints. If $\mathcal S$ is a valid PRS for $U$, then the induced random set
\[
  \Theta_x(\mathcal S)
  =
  \{p:\mathcal S\cap(F_p(x-1),F_p(x)]\neq\emptyset\}
\]
gives calibrated belief and plausibility for assertions about $p$. This preserves Dempster's set-valued insight for discreteness, but it locates the IM prediction problem on the scalar auxiliary scale tied directly to the observed count.

\subsubsection{Bivariate normal correlation}

The bivariate normal correlation coefficient is the cleanest illustration of the local procedure because the usual sufficiency reduction leaves one extra auxiliary dimension. With the notation of Section~\ref{sec:fiducial-confidence}, the reduced association is
\[
  S_1=(1+\rho)U_1,\qquad S_2=(1-\rho)U_2,
  \qquad U_1,U_2\stackrel{\mathrm{ind}}{\sim}\chi^2_n,
\]
where $-1<\rho<1$. Write $s=(s_1,s_2)$ for the observed sufficient statistic. For a fixed candidate value $\rho$, the auxiliary value required by the observation is
\[
  u_{s,\rho}=\left(\frac{s_1}{1+\rho},\frac{s_2}{1-\rho}\right).
\]
The parameter-motion vector field induced by varying $\rho$ while keeping $s$ fixed is therefore
\[
  G^\rho(u)
  =
  \frac{\partial u_{s,\rho}}{\partial\rho}
  =
  \left(-\frac{u_1}{1+\rho},\frac{u_2}{1-\rho}\right)'.
\]
Following \eqref{eq:local-cim-pde}, a local conditioning coordinate $\eta_\rho$ must satisfy
\begin{equation}\label{eq:bvn-corr-local-pde}
  -\frac{u_1}{1+\rho}\frac{\partial\eta_\rho}{\partial u_1}(u)
  +\frac{u_2}{1-\rho}\frac{\partial\eta_\rho}{\partial u_2}(u)
  =0 .
\end{equation}
A corresponding first integral is
\[
  \eta_\rho(u)=(1+\rho)\log u_1+(1-\rho)\log u_2,
\]
since substituting this derivative in \eqref{eq:bvn-corr-local-pde} gives $-1+1=0$. Hence the observed local conditioning value is
\begin{align*}
  h_\rho(s)
  &=\eta_\rho(u_{s,\rho})\\
  &=\log(s_1s_2)+\rho\log(s_1/s_2)
    -\{(1+\rho)\log(1+\rho)+(1-\rho)\log(1-\rho)\}.
\end{align*}
This is the promised local version of $H_{\theta_0}(X)=\eta_{\theta_0}(U)$: the conditioning coordinate is fixed only after the candidate value has been fixed.

The complementary coordinate can be chosen as
\[
  \tau(U)=\log(U_1/U_2).
\]
Then
\[
  T(S)=\log(S_1/S_2)
  =
  z(\rho)+V,\qquad
  z(\rho)=\log\{(1+\rho)/(1-\rho)\},\quad
  V=\tau(U).
\]
Thus, to assess the singleton assertion $\{\rho\}$, the local conditional IM predicts $V$ under the conditional law
\[
  V=\log(U_1/U_2)\mid \eta_\rho(U)=h_\rho(s).
\]
This is exactly the A-step of the local construction: condition on the observed auxiliary feature and retain only the scalar auxiliary variable still needed to check compatibility with the candidate $\rho$.

The conditional density follows by an explicit change of variables. Let
\[
  v=\log(u_1/u_2),\qquad
  h=(1+\rho)\log u_1+(1-\rho)\log u_2 .
\]
Solving gives
\[
  \log u_1=\frac{h+(1-\rho)v}{2},\qquad
  \log u_2=\frac{h-(1+\rho)v}{2}.
\]
Since $U_1$ and $U_2$ are independent $\chi^2_n$ variables, their joint density is proportional to
\[
  (u_1u_2)^{n/2-1}\exp\{-(u_1+u_2)/2\}.
\]
The Jacobian for $(v,h)\mapsto(u_1,u_2)$ is $(u_1u_2)/2$, and
\[
  u_1u_2=\exp(h-\rho v),\qquad
  \frac{u_1+u_2}{2}
  =
  \cosh(v/2)\exp\{(h-\rho v)/2\}.
\]
Therefore, conditional on $h=h_\rho(s)$, the density of $V$ is
\begin{equation}\label{eq:bvn-corr-cond-density}
  f_{\rho,h_\rho}(v)
  \propto
  \exp\left\{-\frac{n\rho v}{2}
  -\cosh(v/2)\exp\left(\frac{h_\rho(s)-\rho v}{2}\right)\right\},
  \qquad v\in\mathbb R.
\end{equation}
Let $F_{\rho,h_\rho}$ denote the distribution function corresponding to \eqref{eq:bvn-corr-cond-density}. For the observed $s$, the auxiliary value required by the candidate $\rho$ in the conditional association is
\[
  c_\rho(s)=T(s)-z(\rho)
  =
  \log(s_1/s_2)-\log\{(1+\rho)/(1-\rho)\}.
\]
Equivalently, $F_{\rho,h_\rho}\{c_\rho(s)\}$ is the scalar uniform auxiliary value associated with the candidate $\rho$ after conditioning. Using the usual symmetric nested PRS on this scalar uniform scale gives the singleton plausibility contour
\[
  \pl_s(\rho)=1-\left|2F_{\rho,h_\rho}\{c_\rho(s)\}-1\right|,
  \qquad -1<\rho<1.
\]

Geometrically, this derivation also shows why the example is not regular. A regular conditioning coordinate would have to solve \eqref{eq:bvn-corr-local-pde} for all $\rho$ simultaneously. But the vector fields $G^\rho(u)$, as $\rho$ varies, span the full two-dimensional tangent space on the positive quadrant, so any common first integral is locally constant. Localization is therefore not a cosmetic choice; it is the operation that leaves a one-dimensional auxiliary prediction problem while preserving conditional validity at the asserted value.

\begin{figure}[tbp]
  \centering
  \includegraphics[width=\textwidth]{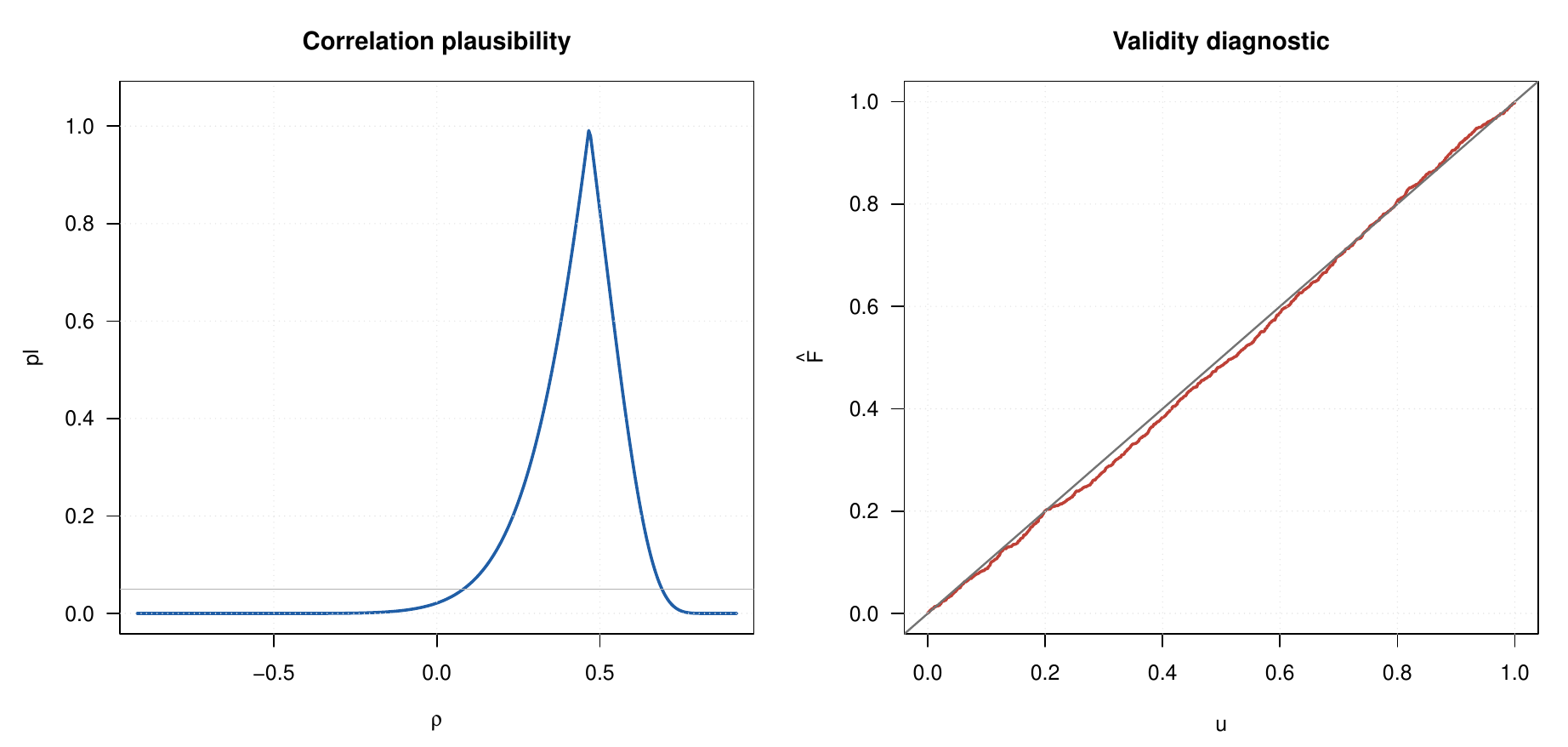}
  \caption{Bivariate normal correlation example. The left panel shows a local conditional plausibility curve for $\rho$ when $n=20$, $s_1=30.60$, and $s_2=10.52$. The right panel is a validity diagnostic based on 1200 simulations under $\rho=0.45$; the empirical distribution of $\pl_S(0.45)$ is close to the uniform benchmark.}
  \label{fig:bvn-correlation-plausibility-validity}
\end{figure}

For a simulated data set with $n=20$ generated at $\rho=0.45$, the observed values are $s_1=30.60$ and $s_2=10.52$. Figure~\ref{fig:bvn-correlation-plausibility-validity} shows the resulting local conditional plausibility curve and a validity diagnostic. Martin and Liu's simulations show that local conditional IM intervals attain the target coverage and can compare favorably with highly accurate conditional likelihood approximations based on $r^\star$ \citep{MartinLiu2015Conditional,BarndorffNielsen1986SignedLogLikelihood,Fraser1990ObservedLikelihood,Reid1995Conditioning}.

\section{Marginal IMs}
\label{sec:marginal}

Scientific questions are usually marginal questions. Researchers ask about a treatment effect, a success probability, a variance component, a risk difference, or a prediction target, not about every unknown in a model simultaneously. Marginal IMs address inference on an interest parameter in the presence of nuisance parameters. The discussion below focuses on \citet{MartinLiu2015Marginal}.

\subsection{Regular Marginal IMs}
\label{sec:mim-methodology}

Let $\theta=(\psi,\xi)$, where $\psi$ is the interest parameter and $\xi$ is nuisance structure. A brute-force IM for the full $\theta$ may be valid, but it often predicts more auxiliary variation than the scientific question requires. This is inefficient, because predictive random sets become larger as the auxiliary dimension grows. The aim of a marginal IM is to construct a lower-dimensional association for $\psi$ itself, then predict only the auxiliary component that remains genuinely uncertain for the targeted inference.

The regular case of \citet{MartinLiu2015Marginal} makes the idea precise. Suppose a baseline association can be rewritten, after a smooth transformation of data and auxiliary variables, in the form
\[
  \bar p(X,\psi)=\bar a(V_1,\psi), \qquad
  c(X,V_2,\psi,\xi)=0,
\]
where $V=(V_1,V_2)$ has known distribution. The first equation involves the interest parameter and the auxiliary feature $V_1$; the second equation contains nuisance structure. The association is regular, in the sense of \citet{MartinLiu2015Marginal}, when the nuisance equation can be solved in $\xi$ for any fixed $(x,v_2,\psi)$ in the relevant range. Then $V_2$ carries no direct information about $\psi$. It can be ignored, or equivalently the predictive random set can be stretched to fill the entire $V_2$ direction, without damaging validity for $\psi$.

The resulting marginal association for inference about $\psi$ is
\[
  \bar p(X,\psi)=\bar a(V_1,\psi), \qquad V_1\sim P_{V_1}.
\]
For observed $x$ and auxiliary value $v_1$, define
\[
  \Psi_x(v_1)=\{\psi:\bar p(x,\psi)=\bar a(v_1,\psi)\}.
\]
If $\mathcal S$ is a valid predictive random set for $V_1$, then
\[
  \Psi_x(\mathcal S)=\bigcup_{v_1\in\mathcal S}\Psi_x(v_1)
\]
induces marginal belief and plausibility functions
\[
  \operatorname{mbel}_x(A)=P_{\mathcal S}\{\Psi_x(\mathcal S)\subseteq A\},
  \qquad
  \operatorname{mpl}_x(A)=P_{\mathcal S}\{\Psi_x(\mathcal S)\cap A\neq\emptyset\}.
\]
The validity theorem says, roughly, that valid prediction of $V_1$ gives valid marginal inference for $\psi$, uniformly over nuisance values $\xi$. This is the key point: nuisance elimination happens before probability calculus is transferred to the parameter space.

Like inference with basic IMs, it is often useful to summarize a marginal IM by a pointwise plausibility contour
\[
  \psi\mapsto \pl_x(\psi).
\]
Its level sets $\{\psi:\pl_x(\psi)>\alpha\}$ give plausibility regions when the pointwise construction has the appropriate validity property. For more general assertions $A\subseteq\Psi$, a consonant or possibilistic extension is natural:
\[
  \pl_x(A)=\sup_{\psi\in A}\pl_x(\psi),\qquad
  \bel_x(A)=1-\sup_{\psi\notin A}\pl_x(\psi).
\]

\subsection{Non-regular marginal IMs}
Regularity is sufficient but not universal. In nonregular problems, the reduced auxiliary component relevant to $\psi$ may still have a distribution depending on $\xi$. Then simply ignoring nuisance structure can lose validity. The marginal IM response is to use a uniformly valid predictive random set: a random set that predicts the nuisance-dependent auxiliary variable validly for every possible nuisance value. One practical route is to replace the nuisance-dependent auxiliary variable by a nuisance-free stochastic bound with fatter tails, then predict that bound. The resulting generalized marginal IM may be conservative, but the conservatism is explicit and protects the plausibility scale. These are illustrated in detail below using the Behrens-Fisher example.

A differential geometric formulation of non-regular marginal IMs, including regular marginal IMs as special cases, provides an elegant mathematical theory. The details are given in \citet{Liu2026DifferentialGeometricIMs}.

\subsection{Examples}
\label{sec:mim-examples}

\subsubsection{A Regular Marginal IM Example}
\label{sec:regular-mim-examples}

The unknown-variance normal example is the cleanest illustration. If $\psi=\mu$ and $\xi=\sigma$, the conditional reduction to $(\bar X,S)$ gives
\[
  \bar X=\mu+\sigma n^{-1/2}U_1,\qquad S=\sigma U_2,
\]
where $U_1\sim N(0,1)$ and $(n-1)U_2^2\sim\chi^2_{n-1}$, independently. Dividing the centered mean equation by the observed scale removes $\sigma$ and yields the marginal association
\[
  \frac{\bar X-\psi}{S/\sqrt n}=T,\qquad T\sim t_{n-1}.
\]
The plausibility contour for $\psi$ follows directly from this one-dimensional auxiliary variable. This is what a good marginal procedure should do: answer the question about $\mu$ without requiring a prior distribution for $\sigma$ and without first constructing a joint posterior-like object for $(\mu,\sigma)$. The success here should not be overgeneralized. It works because the nuisance scale can be removed by a pivot whose distribution is free of $\sigma$.

\subsubsection{A Nonregular Marginal IM Example}
\label{sec:behrens-fisher-mim}

The Behrens--Fisher problem is the classical two-sample normal problem with
unknown and unequal variances, originating in the work of
\citet{Behrens1929FewObservations} and Fisher's fiducial analysis
\citep{Fisher1935Fiducial}.  It has since served as a useful testing ground
for fiducial, Bayesian, frequentist, and more recently IM approaches; see,
for example, the practical discussion in \citet{Scheffe1970BehrensFisher},
the review of \citet{KimCohen1998BehrensFisherReview}, the exact-solution
comparison of \citet{DudewiczMaMaiSu2007BehrensFisher}, and the broader
review by \citet{PaulWangUllah2019BehrensFisherReview}.  The marginal IM
treatment in \citet{MartinLiu2015Marginal} is especially relevant here
because it shows how the example exposes the boundary between regular and
nonregular marginalization. What is discussed below follows \citet{WangLiu2026BehrensFisher}. Let
\[
  X_{ki}\stackrel{\mathrm{ind}}{\sim}N(\mu_k,\sigma_k^2),
  \qquad i=1,\ldots,n_k,\quad k=1,2,
\]
with interest focused on
\[
  \psi=\mu_2-\mu_1 .
\]
The nuisance structure consists of the remaining location direction and the two
variances.  This familiar example is useful here because regular
marginalization almost works, but not quite.

After the conditional-IM reduction to sample means and sample standard
deviations, the retained association is
\[
  \bar X_k=\mu_k+\sigma_k n_k^{-1/2}U_{1k},
  \qquad
  S_k=\sigma_k U_{2k},\qquad k=1,2,
\]
where $U_{1k}\sim N(0,1)$ and $(n_k-1)U_{2k}^2\sim\chi^2_{n_k-1}$, independently.
The reduction from the original four parameters to the two-parameter
marginalization association proceeds in three steps. First, reparameterize
\[
  \lambda=\mu_1,\qquad
  \psi=\mu_2-\mu_1,\qquad
  \omega=f(\sigma_1,\sigma_2),\qquad
  \xi=\frac{\sigma_1^2/n_1}{\sigma_1^2/n_1+\sigma_2^2/n_2},
\]
where
\[
  f(a,b)=\left(\frac{a^2}{n_1}+\frac{b^2}{n_2}\right)^{1/2},
\]
so that
\[
  \mu_1=\lambda,\qquad \mu_2=\lambda+\psi,\qquad
  \sigma_1=\omega(n_1\xi)^{1/2},\qquad
  \sigma_2=\omega\{n_2(1-\xi)\}^{1/2}.
\]
Here $\lambda$ is the irrelevant location level, $\omega$ is the overall
standard-error scale, and $\xi\in(0,1)$ is the nuisance variance-ratio
coordinate.

Second, subtract the two sample-mean associations. With
$\bar Y=\bar X_2-\bar X_1$,
\[
  \bar Y
  =
  \psi
  +\omega\{\sqrt{1-\xi}\,U_{12}-\sqrt{\xi}\,U_{11}\}.
\]
For each fixed $\xi$, the auxiliary combination
\[
  U_1(\xi)=\sqrt{1-\xi}\,U_{12}-\sqrt{\xi}\,U_{11}
\]
has a $N(0,1)$ distribution and is independent of $(U_{21},U_{22})$.
Thus the nuisance location $\lambda$ has disappeared, but the overall scale
$\omega$ and the variance-ratio coordinate $\xi$ remain.

Third, remove $\omega$ by studentizing with the observed standard errors. Since
\[
  \frac{S_1^2}{n_1}+\frac{S_2^2}{n_2}
  =
  \omega^2\{\xi U_{21}^2+(1-\xi)U_{22}^2\},
\]
one obtains
\[
  T_\psi(X)=
  \frac{\bar X_2-\bar X_1-\psi}
       {\{S_1^2/n_1+S_2^2/n_2\}^{1/2}} .
\]
This statistic satisfies
\[
  T_\psi(X)
  =
  \frac{U_1(\xi)}
       {\{\xi U_{21}^2+(1-\xi)U_{22}^2\}^{1/2}} .
\]
The remaining information about $\xi$ comes from the observed ratio of sample
standard errors:
\[
  R=\frac{n_1S_2^2}{n_2S_1^2}
  =
  \frac{1-\xi}{\xi}\frac{U_{22}^2}{U_{21}^2}.
\]
Writing $U_1$ for a generic $N(0,1)$ variable with the same distribution as
$U_1(\xi)$, the resulting two-parameter marginalization association for
$(\psi,\xi)$ is
\begin{align}
  T_\psi(X)
  &=
  Z_1(\xi)
  =
  \frac{U_1}{\{\xi U_{21}^2+(1-\xi)U_{22}^2\}^{1/2}},
  \label{eq:behrens-fisher-z1}\\
  R
  &=
  Z_2\,\frac{1-\xi}{\xi},
  \qquad
  Z_2=\frac{U_{22}^2}{U_{21}^2}.
  \label{eq:behrens-fisher-z2}
\end{align}
The location $\lambda$ and overall scale $\omega$ have been eliminated before
any predictive random set is introduced. If the law of $Z_1(\xi)$ were free of
$\xi$, this would be a regular marginal association for $\psi$. Instead, the
law of $Z_1(\xi)$ changes with $\xi$: at the variance-ratio boundaries it
approaches Student $t_{n_1-1}$ or $t_{n_2-1}$ limits, and in the balanced case
the interior law is still different from the common boundary law.  The second
equation does not eliminate the nuisance ratio; solving it gives
$\xi=Z_2/(Z_2+R)$, which depends on the unobserved auxiliary $Z_2$.  Thus the
exact reduced association remains genuinely two-dimensional.

Following the validity-first treatment in \citet{WangLiu2026BehrensFisher},
which builds on the generalized marginal IM strategy of
\citet{MartinLiu2015Marginal}, this nonregularity is handled by predicting a
two-dimensional auxiliary variable, then using a projection that is sharp for
inference on the marginal contrast.  The key probability
bound is due to Hsu \citep{Hsu1938TTest} and was later emphasized as a
practical conservative solution by \citet{Scheffe1970BehrensFisher}.  It says
that, for
\[
  m=\min(n_1,n_2)-1,
\]
the symmetric tails of $Z_1(\xi)$ are bounded uniformly by those of a Student
$T_m$ variable:
\[
  \pr_\xi\{|Z_1(\xi)|>c\}
  \leq
  \pr\{|T_m|>c\},
  \qquad c\geq0.
\]
The bound is sharp at the least favorable variance-ratio boundary.  The
corresponding IM predictive random set is the cylinder
\begin{equation}
  \mathcal S_{\mathrm{IM}}=[-C,C]\times\mathbb R_+,
  \qquad C\sim |T_m|.
  \label{eq:behrens-fisher-cylinder}
\end{equation}
It is sharp in the mean-contrast coordinate and vacuous in the variance-ratio
coordinate.  The induced plausibility contour is
\begin{equation}\label{eq:behrens-fisher-plausibility}
  \mpl_x(\psi)
  =
  2\left[
  1-G_m\left(
    \left|
    \frac{\bar x_2-\bar x_1-\psi}
         {\{s_1^2/n_1+s_2^2/n_2\}^{1/2}}
    \right|
  \right)
  \right],
\end{equation}
where $G_m$ is the cdf of $T_m$.  The associated
$100(1-\alpha)\%$ plausibility interval is the Hsu--Scheffe interval
\[
  \bar x_2-\bar x_1
  \pm
  t_{m,1-\alpha/2}
  \left(\frac{s_1^2}{n_1}+\frac{s_2^2}{n_2}\right)^{1/2},
\]
and it is valid uniformly over the unknown variances.

The important lesson is not only validity but the order in which validity and
efficiency are judged.  Since the nuisance ratio is not removed by regular
marginalization, any prior-free fixed-sample method must either use a uniform
bound or introduce additional structure.  In the cylindrical class, Hsu's sharp
boundary bound makes the IM interval minimax and admissible: no smaller
two-sided critical value can retain uniform validity.  More generally, marginal
inference on $\psi$ depends only on the first-coordinate projection of a
two-dimensional PRS, so attempts to predict the variance-ratio auxiliary can
redistribute interval width across variance-ratio regimes but cannot uniformly
shorten the projection without paying elsewhere.

This explains the relation with classical competitors.  Welch's approximation
\citep{Welch1938UnequalVariances,Welch1947GeneralizationStudent}, together
with Satterthwaite's moment-matched degrees of freedom
\citep{Satterthwaite1946VarianceComponents}, gives the familiar
Welch--Satterthwaite interval.  It uses the observed standard errors to choose
a data-dependent reference $t$ distribution and is often shorter, but it is not
an exact finite-sample procedure with uniform coverage over the variance ratio.
Fisher's fiducial, generalized fiducial, and independent-Jeffreys Bayesian
intervals use the observed ratio of standard errors in different ways; they can
be reasonable, and in special comparisons even attractive, but their efficiency
comparisons depend on the criterion used to weight variance-ratio regimes.  The
IM/Hsu answer is the validity-first answer: it is conservative because it
refuses to average over or otherwise guess the nuisance variance ratio.  The
example therefore shows both the value and the limit of marginal IMs.  Auxiliary
variables can diagnose exactly where nuisance information remains, and validity
can be preserved without a prior, but some nonregular problems force
conservatism if the guarantee is to be uniform.

\section{Conclusion}
\label{sec:conclusion}

The main message is that scientific uncertainty should be both data-conditional and frequency-calibrated. The auxiliary-variable perspective connects these two demands when a genuine prior is not available. Inferential models continue a long historical project on inference with auxiliaries. Fisher developed fiducial inference as an attempt to obtain post-data probability statements about unknown parameters directly from the sampling model, by transferring the known randomness of pivotal or auxiliary quantities to the parameter space without introducing a prior distribution. Neyman supplied mathematical clarity through confidence procedures. Dempster and Shafer supplied a language of belief, plausibility, and lower and upper probabilities. Generalized fiducial inference revived the auxiliary-variable mechanism. IMs synthesize these ideas by treating the unobserved auxiliary variable as the carrier of inferential uncertainty and predicting it in a calibrated way, with efficiency pursued under the constraint of validity.

Several extensions and applications deserve further reading. For imprecise-probability updating and the IM response to unsettling combination rules, see \citet{GongMeng2021Unsettling,LiuMartin2021SettleUnsettling}; for p-values as plausibilities, see \citet{MartinLiu2014Pvalues}; for the confidence-curve and fiducial reconstruction viewpoint, see \citet{CuiHannig2025Demystifying,XieSingh2013ConfidenceDistribution,SchwederHjort2016}.  For modern generalized and possibilistic IMs based on calibrated contours, see \citet{LiuMartin2020Possibility,Martin2026PossibilisticReview,Martin2015PlausibilityExact,Martin2018GeneralizedAssociations,CahoonMartin2021Censored,JiangLiuZhang2024CalibratedBootstrap,Martin2025MonteCarloPossibilistic}; related betting-IM and e-process constructions are developed in \citet{Martin2024WaudbyRamdasDiscussion,Martin2026RegularizedEProcesses,DixitMartin2025PredictiveRecursion,DeyMartinWilliams2025GeneralizedUniversal}. For non-parametric and high-dimensional directions, many-normal-means and shape-constrained examples provide useful entry points, including high-dimensional structural learning \citep{EschkerLiu2026RevisitingStein} and partial conditioning for H{\"o}lder-constrained means \citep{YangWangLiu2023PartialConditioning}.

For readers shaped by Jaynes' view of probability as the logic of science, IMs should not be read as a rejection of that logic \citep{Jaynes2003LogicScience}. They sharpen its starting point. When the sampling model and a genuine prior distribution are both on the table, Bayes' theorem is the natural way to refine knowledge through the observed data. The Martin--Liu question is what the logic of science should do when the prior entry is empty. The IM answer is not automatically to fill that empty place with a formal or default prior and then interpret the resulting posterior as objective scientific uncertainty. It is to keep on the table what is still genuinely available: the sampling model, the auxiliary variable, its known distribution, the observed association, and the scientific assertion, some of which become redundant and, thereby, often hidden when a prior is available. In this sense, Jaynes teaches that probability is the right logic when the information warrants a probability distribution; Martin and Liu teach that the logic of science must also say what to do when it does not. Thus, the development of IMs using auxiliary variables can be read as continuing Jaynes' project of probability as the logic of science, but with an important qualification: when the available information does not justify a full prior probability distribution, the logic of uncertainty should be calibrated, auxiliary-variable based, and sometimes imprecise. 

There is still much to do beyond expanding the range of IM applications. For continuous-data models, a differential-geometric theory of IMs now seems within reach: conditioning, marginalization, auxiliary reductions, and local validity all have natural geometric interpretations in terms of transformations, fibers, transverse coordinates, and induced contours. Discrete-data models are more challenging because the clean geometry of smooth auxiliary maps is replaced by jumps, nonunique inversions, and boundary-sensitive validity. More broadly, since IMs use auxiliary variables as principled surrogates when prior distributions are unavailable, a new look at the likelihood principle may be possible. Such a reconstruction could clarify why likelihood alone is often compelling but not sufficient for calibrated uncertainty, and it may shed light on long-standing debates concerning the likelihood principle and the foundations of statistical inference more generally.

\section*{Acknowledgments}

This work was supported in part by the U.S. National Science Foundation under grant DMS-2412629.

\bibliographystyle{chicago}
\bibliography{references}

\end{document}